\documentclass{amsart}
\usepackage{amsmath, amsthm, amssymb, graphicx, tikz, zref} 
\usetikzlibrary{arrows}
\usepackage{enumitem}
\usepackage{courier}

\theoremstyle{definition}
\newtheorem{remark}[equation]{Remark}

\newtheorem{theorem}{Theorem}[section]

\newtheorem{conjecture}[theorem]{Conjecture}
\newtheorem{defn}[theorem]{Definition}

\usepackage[margin=1in]{geometry} 
\usepackage{tikz-cd}

\newcommand{\RA}{\textcolor{magenta}}

\usepackage[foot]{amsaddr}

\newcommand{\Hom}{\mathrm{Hom}}
\usepackage{placeins}
\DeclareUnicodeCharacter{2212}{-}
\usepackage{comment}

\usepackage{pdfpages}
\usepackage{hyperref}
\usepackage{faktor}

\usepackage[backend=bibtex, 
style=numeric, natbib=true
]{biblatex}
\usepackage{tikz-cd}
\usepackage{mathabx,epsfig}

\hypersetup{
    colorlinks,
    linkcolor={red!50!black},
    citecolor={blue!50!black},
    urlcolor={blue!80!black}
}
\newcommand{\R}{\mathbb{R}}

\newcommand{\PP}{\mathbb{P}}\newcommand{\Z}{\mathbb{Z}}\newcommand{\C}{\mathbb{C}}

\usepackage{makecell}

\title{Exceptional Collections for Toric Fano Fivefolds}
\author{Reginald Anderson}
\address[Reginald Anderson]{UC Irvine}

\email{reginala@uci.edu }

\date{July 2025}

\begin{document}

\maketitle

\begin{abstract}
Resolutions of the diagonal of toric varieties has been an active area of study since Beilinson's celebrated resolution of the diagonal for $\PP^n$ and the disproof of King's conjecture. The author generalized a cellular resolution of the diagonal given by Bayer-Popescu-Sturmfels to yield a virtual resolution of the diagonal for smooth projective toric varieties, which extends to toric Deligne-Mumford stacks which are a global quotient of a smooth projective variety by a finite abelian group. Moreover, a celebrated result of Hanlon-Hicks-Lazarev gives a symmetric, minimal resolution of the diagonal for smooth projective toric varieties. This work studies when smooth projective toric Fano varieties in dimension 5 yield exceptional collections of line bundles using resolution of the diagonal. We give the first known count of 300 out of 866 smooth projective toric Fano 5-folds for which the Hanlon-Hicks-Lazarev resolution of the diagonal yields a full strong exceptional collection of line bundles. 

\end{abstract}
\vspace{.2cm} 

\section{Introduction}

 The gold standard for a presentation of the bounded derived category of coherent sheaves on a variety $X$, $\mathrm{D}^b_{Coh}(X)$, is a \emph{full strong exceptional collection} \cite{Baer88, BondalRepnAssocAlg}. 
The utility of having such a decomposition was first observed when Beilinson gave a resolution of the diagonal for $\mathbb{P}^n$; concretely describing $\mathrm{D}^b_{Coh}(\mathbb{P}^n)$ using subcategories generated by line bundles on $X$ vastly simplifies homological computations (K-theory, stability conditions, etc.). In general, full strong exceptional collections are elusive. Even if they exist, as Kawamata showed for smooth projective toric varieties,  the possible structure of exceptional objects is widely unknown \cite{Kawamata2017}. King conjectured \cite{KingUnpub} in unpublished notes that any smooth projective toric variety has a full strong exceptional collection of line bundles. King's conjecture was proven false by Hille-Perling and Efimov \cite{Efimov2014, hille2006counterexample}. Presently, whether or not a given smooth projective toric variety admits a full strong exceptional collection of line bundles is checked on a case-by-case basis; the failure of one method of detecting the existice of a full strong exceptional collection of line bundles does not imply a full strong exceptional collection of line bundles does not exist. 
Full strong exceptional collections of line bundles on all 124 smooth projective toric Fano fourfolds were constructed by Prabhu-Naik in \cite{PRABHUNAIK2017348} and used to create tilting bundles, building on work of Craw \cite{craw2008quiverrepresentationstoricgeometry}.

For a smooth projective Fano or variety of general type, $\mathrm{D}^b_{Coh}(X)$ determines $X$ up to isomorphism \cite{bondal2002derivedcategoriescoherentsheaves}. Beilinson's resolution of the diagonal \cite{Beilinson1978} gave a full strong exceptional collection of line bundles for $\PP^n(\C)$, which is a smooth projective toric Fano variety. Generalizing Beilinson's resolution of the diagonal, Bayer-Popescu-Sturmfels \cite{bayer-popescu-sturmfels} gave a cellular resolution \cite[Ch. 4]{miller-sturmfels} of the diagonal for a proper subclass of smooth projective toric varieties which they called ``unimodular," which is a more restrictive condition than being nonsingular. Generalizing this, a non-minimal virtual \cite[Definition 1.1]{Erman_2020} cellular resolution of the diagonal was given by the author in \cite{anderson2023resolutiondiagonalsmoothtoric} for smooth projective toric varieties and in \cite{anderson2023resolutiondiagonaltoricdelignemumford} for toric Deligne-Mumford stacks which are a global quotient of a smooth projective toric variety by a finite abelian group. In \cite{hanlon2023resolutions}, Hanlon-Hicks-Lazarev gave a minimal virtual resolution of the diagonal for any toric subvariety of a smooth projective toric variety. 


It is proven that the collection of line bundles appearing from the Hanlon-Hicks-Lazarev resolution of the diagonal agrees with the Bondal-Thomsen collection for smooth projective toric Fano varieties. In \cite{bondalOberwolfach}, Bondal gives a numerical criterion which holds for $16/18$ smooth projective toric Fano threefolds. This has been checked by the author to coincide with the success of the Hanlon-Hicks-Lazarev resolution to yield a full strong exceptional collection of line bundles in dimension 3, see  \cite{anderson2024exceptionalcollectionslinebundles}. There are also computational results in dimension 4 \cite{ramirez2025exceptionalcollectionstoricfano}. Smooth projective toric varieties for which the Hanlon-Hicks-Lazarev resolution yields a full strong exceptional collection of line bundles are referred to as Bondal-Ruan type \cite{favero2022homotopy} \cite{bondalOberwolfach}. This raises the question: ``What happens in higher dimension?" As a first step towards this question, we conclude that for dimension 5:

\begin{theorem} \label{thm: Thm1} The Hanlon-Hicks-Lazarev resolution of the diagonal yields a full strong exceptional collection of line bundles for 300 out of 866 smooth projective toric Fano fivefolds.  \end{theorem}

To the author's knowledge, this is the first work giving the count of 300 out of 866 smooth projective toric Fano fivefolds admitting a Hanlon-Hicks-Lazarev resolution of the diagonal. This is in contrast to dimension 4, where multiple methods of determining the count of 72/124 fourfolds whose Hanlon-Hicks-Lazarev resolution of the diagonal imply the existence of full strong exceptional collections of line bundles have been used \cite{4foldsundergrads,PRABHUNAIK2017348}. In general, the count of smooth toric Fano $n$-folds admitting a Hanlon-Hicks-Lazarev resolution of the diagonal were previously known for $n\le 4$. 

Smooth projective toric varieties for which the Hanlon-Hicks-Lazarev resolution yields a full strong exceptional collection of line bundles are referred to as Bondal-Ruan type \cite{bondalOberwolfach,favero2022homotopy}. 

The indices q for which the Hanlon-Hicks-Lazarev resolution of the diagonal yields a full strong exceptional collection of line bundles on 
\texttt{
    X = smoothFanoToricVariety(5,q)}
is given by the following list:
\begin{align*} \{&0, 476, 477, 478, 479, 480, 481, 482, 483, 484, 485, 486, 487, 488, 489, 490, 491, 492, 493, 494, 496, 497, 498, \dots \\
\dots & 499, 500, 502, 503, 504, 505, 506, 507, 508, 509, 510, 511, 512, 513, 514, 515, 517, 518, 519, 520, 521, 522, 523, \dots \\
\dots & 524, 525, 526, 527, 528, 529, 530, 531, 532, 533, 534, 536, 537, 538, 539, 540, 541, 542, 543, 544, 545, 546, 547, \dots \\
\dots & 551, 552, 553, 554, 555, 556, 559, 561, 562, 563, 564, 565, 566, 567, 570, 571, 572, 573, 574, 575, 576, 577, 579, \dots \\
\dots & 581, 582, 583, 584, 585, 586, 587, 588, 590, 591, 592, 593, 594, 595, 599, 600, 601, 602, 603, 604, 606, 608, 609, \dots \\
\dots & 610, 611, 612, 613, 614, 615, 617, 620, 621, 623, 624, 625, 626, 628, 629, 630, 631, 632, 633, 634, 637, 638, 639, \dots \\
\dots & 642, 643, 646, 647, 648, 649, 650, 651, 652, 653, 655, 656, 658, 660, 661, 662, 663, 664, 665, 667, 668, 669, 670, \dots \\
\dots & 671, 672, 673, 675, 676, 677, 678, 679, 680, 681, 682, 685, 686, 687, 688, 689, 690, 691, 692, 693, 694, 695, 697, \dots \\
\dots & 698, 699, 700, 702, 703, 705, 706, 709, 710, 712, 713, 714, 715, 716, 717, 718, 719, 720, 721, 722, 724, 725, 730, \dots \\
\dots & 731, 732, 733, 734, 735, 737, 738, 739, 740, 743, 744, 745, 746, 747, 748, 749, 750, 751, 752, 753, 754, 758, 760, \dots \\
\dots & 763, 766, 767, 773, 774, 775, 776, 777, 778, 779, 780, 781, 782, 784, 785, 788, 789, 791, 792, 794, 795, 796, 797, \dots \\
\dots & 799, 800, 801, 802, 805, 809, 810, 811, 812, 813, 814, 816, 817, 818, 819, 820, 821, 822, 823, 824, 825, 826, 829, \dots \\
\dots & 830, 831, 832, 834, 835, 837, 840, 841, 848, 849, 850, 851, 852, 853, 854, 855, 856, 857, 858, 859, 861, 862, 863, 864 \}.
\end{align*}

We make this list available at 

\begin{center}\url{https://github.com/reggiea91/ADJOINT2025_Fivefolds.git}\end{center}

\begin{remark} Though we work exclusively over $\C$, full strong exceptional collections of line bundles on toric varieties have also been considered over more general fields \cite{ballard2020derivedcategoriescentrallysymmetricsmooth} \cite{BallardDuncanMcFaddin}.
\end{remark}

\section{Background}

We work over the algebraically closed field $\Bbbk\cong\overline{\Bbbk}\cong \C$. A smooth projective toric variety $X$ is said to be Fano if the 
anticanonical bundle is ample. Let $\mathcal{C}$ denote a $\Bbbk$-linear triangulated category. Recall \cite{huybrechts2006fourier} that in a $k$-linear triangulated category $\mathcal{D}$, an object $E$ is called \textbf{exceptional} if \[ \mathrm{Hom}_\mathcal{D}(E, E[\ell]) = \begin{cases} \Bbbk & \ell = 0, \\ 0 & \text{ else.} \end{cases} \] 

An \textbf{exceptional sequence} $\mathcal{E}$ is a sequence $E_1, \dots, E_r$ of exceptional objects such that $\Hom_\mathcal{D}(E_i, E_j[\ell]) = 0$ for $i>j$ and all $\ell$. That is, 

\[ \text{Hom}_\mathcal{D}(E_i, E_j[\ell]) = \begin{cases} \Bbbk & i=j, \ell=0 \\ 0 & \text{ if } i>j \text{ or if } \ell\neq 0, i=j. \end{cases} \] 

An exceptional sequence is \textbf{full} if $\mathcal{D}$ is generated by $\{E_i\}$. An exceptional collection $E_1, \dots, E_r$ is \textbf{strong} if \[ \text{Hom}_{\mathcal{D}}(E_i, E_j[\ell]) = 0\] for all $i,j$, and $\ell\neq 0$. Here, $\mathcal{D} = \mathrm{D}^b_{Coh}(X)$, the bounded derived category of coherent sheaves on a smooth projective toric Fano variety of dimension $5$ over $\mathbb{C}$. Given a full strong exceptional collection of line bundles $\mathcal{E} = \{E_i\}_{i=1}^r$ on a smooth projective variety, we form $\mathcal{G} = \bigoplus_{i=1}^r E_i$ and $A =\text{End}_{\mathrm{D}^b_{Coh}(X)}(\mathcal{G})$ so that $\mathrm{D}^b_{Coh}(X) \simeq \mathrm{D}^b(\text{mod}-A)$ \cite{BondalRepnAssocAlg}. 

\section{Methods}
We share code reproducing our results online at 

\begin{center} \label{link: githublink} \url{https://github.com/reggiea91/ADJOINT2025_Fivefolds.git}\end{center}

 Our code does generalize immediately to higher dimensions for which similar databases are available; while the Macaulay2 database  ends at dimension 6 \cite{M2database}, the polymake database goes through dimension 9 \cite{polymakedatabase}. One would need to modify our code for the polymake database to iterate over the polymake database, and modify accordingly to follow their conventions. Implementing our code also becomes an issue of computing power in order to compile and reach a final count in higher dimensions. Our work in dimension 5 was made possible largely by access to the high performance computer at the Simons Laufer Mathematics Institute, and dimension 6 would likely require greater  prolonged access to such a high performance computer. 

The methods that we have used for this project rely on the Macaulay2 database of smooth projective toric Fano varieties of dimension 5. This database is accessible from the NormalToricVarieties package via 
\texttt{ X = smoothFanoToricVariety(5,q)} for a given index q in $[0, 865]$. Given $X$, a smooth projective toric Fano variety of dimension $5$, we apply the Hanlon-Hicks-Lazarev resolution to $\mathcal{O}_\Delta$ for $X \cong \Delta \subset X\times X$ the diagonal embedding. This yields a locally-free resolution of $\mathcal{O}_\Delta$ in $\mathrm{D}^b(X \times X)$ in terms of box products of line bundles. The resulting collection is symmetric in the sense that the collection of line bundles on the left- and right-hand sides, respectively, are equal. We therefore consider the collection of line bundles $\mathcal{E}$ appearing on the left-hand side and ask whether there exists an ordering of the collection of line bundles $\mathcal{E} = \{E_i\}_{i\in I}$ for which $\mathcal{E}$ gives a full exceptional collection, and further, whether this collection is strong.

In order to check for the existence of such an ordering, we construct a directed graph $G$ of all homomorphisms between line bundles in $\mathcal{E}$ and check for the presence of directed cycles (apart from the self-loops at each vertex to itself). That is, we check for the presence of directed cycles of length greater than $1$ in $G$. The absence of directed cycles of length greater than $1$ in $G$ then implies the existence of an ordering for which the collection $\mathcal{E}$ is exceptional. 

Here, $G$ shows $ \Hom_{\mathrm{D}^b_{Coh}(X)}(E_i, E_j[\ell])$ 
for all $i,j$, and $\ell$. To compute $\Hom_{\mathrm{D}^b_{Coh}(X)}(E_i, E_j[\ell])$, we make use of the fact that for line bundles $\mathcal{O}(D)$ and $\mathcal{O}(E)$ on a smooth projective toric variety $X$, 

\[ \Hom_{\mathrm{D}^b_{Coh}(X)}(\mathcal{O}(D), \mathcal{O}(E)) \cong H^*(X, \mathcal{O}(E-D) ). \]

Bondal's numerical criterion checks for whether the tilting object $\mathcal{G} = \bigoplus_{i\in I} E_i$ (without repetitions) has any higher self-Ext's, so use the presence of an ordering for which $\mathcal{E}$ is exceptional to also indicate that $\mathcal{E}$ is strong. That is, we do not speak of the collection $\mathcal{E}$ being strong if $\mathcal{E}$ is not also exceptional. The rank of graded homomorphisms between the objects $E_i$ in $\mathcal{E}$ is given in an adjacency matrix for $G$ as a lower-triangular matrix whose $(i,j)$ entry indicates the rank of $\Hom_{\mathrm{D}^b_{Coh}(X)}(E_j, E_i[0])$. 

Examples are given for the positive example 

\begin{center}
    \texttt{X = smoothFanoToricVariety(5,851)}
\end{center}

and for the negative example 

\begin{center}
    \texttt{X = smoothFanoToricVariety(5,200)}\end{center}



The author constructed a virtual cellular resolution of the diagonal for smooth projective toric varieties in \cite{anderson2023resolutiondiagonalsmoothtoric} 
generalizing Bayer-Popescu-Sturmfels' approach \cite{bayer-popescu-sturmfels}, which in turn generalized Beilinson's resolution of the diagonal for $\PP^n$ \cite{Beilinson1978}. Smooth projective toric varieties have a presentation of the class group using integer-valued matrices via the fundamental exact sequence \[ 0 \rightarrow M \stackrel{B}{\rightarrow} \Z^{|\mathcal{A}|} \stackrel{\pi}{\rightarrow} \text{Cl}(X) \rightarrow 0 \] where $M$ is the character lattice of the Zariski-dense open algebraic torus $T \cong (\C^*)^m \subseteq X$ of rank $m$, $\mathcal{A}=\Sigma(1)$ is the set of primitive ray generators with $|\Sigma(1)|=n$, and, fixing a basis for $M$, $B$ is represented by the $n \times m$ matrix with rows as primitive ray generators, which are the elements of $\mathcal{A}$. Let $L$ denote the image of $B$ in $\Z^{|\mathcal{A}|}$, and consider $\R L = L\otimes_{\Z} \R \subset \R^{|\mathcal{A}|} = \Z^{|\mathcal{A}|} \otimes_\Z \R$. $\R L$ and the infinite hyperplane arrangement $\mathcal{H} = \{ x_i = j \text{ }|\text{ } 1 \leq i \leq n, j \in \Z\} \subset \R^{|\mathcal{A}|}$ were used in \cite{anderson2023resolutiondiagonalsmoothtoric} to construct an infinite and quotient cellular complex, respectively, to build a cellular resolution of the diagonal. There, the Laurent monomial labeling on the vertex $p= (p_1, \dots, p_n)$ in the finite quotient cellular complex $\faktor{\mathcal{H}_L}{L}=\faktor{\R L \cap \mathcal{H}}{L}$ is \[ \frac{x^{\lfloor p \rfloor} }{y^{\lfloor p \rfloor}}. \]

Hanlon-Hicks-Lazarev \cite{hanlon2023resolutions} use the same hyperplane arrangement and a different convention on Laurent-monomial labelings for the finite quotient cellular complex which agrees with \cite{bondalOberwolfach} to build a symmetric minimal resolution of the diagonal \cite{brown2023short} for smooth projective toric varieties as follows. Given a toric morphism $f: X \rightarrow Y$ of smooth projective toric varieties, Hanlon-Hicks-Lazarev construct a resolution of $f_* \mathcal{O}_X$ considered as an object of $\mathrm{D}^b(Y)$. In the case of the diagonal morphism $\Delta: X \rightarrow X \times X$, $\Delta_* \mathcal{O}_X \cong \mathcal{O}_\Delta$. The diagonal morphism of toric varieties induces a dual map of (stacky) fans $f^\vee: \Sigma \times \Sigma \rightarrow \Sigma$ with kernel of the dual map of lattices isomorphic to the antidiagonal \[ \{(u,-u) \text{ }|\text{ }u \in L \} \] for $L$ the same image of $B$ in the notation of Bayer-Popescu-Sturmfels above. This kernel of the dual map of lattices is isomorphic to the same antidiagonal lattice which is the Lawrence lifting of the lattice $L$ in \cite{bayer-popescu-sturmfels}. Bondal's perspective on toric mirror symmetry \cite{bondalOberwolfach} views the above cellular complex $\faktor{\mathcal{H}_L}{L}$ as the stratification of a real torus $\mathbb{T}$, and gives a morphism $\Phi: \mathbb{T} \rightarrow \mathrm{Pic}(X)$ by \[ \Phi(\sum a_i D_i) = -\sum \lfloor a_i \rfloor D_i \]

For the antidiagonal lattice, Hanlon-Hicks-Lazarev (in their notation) accordingly prescribe the monomial labeling to any $m\in M_\R$

\[ m \mapsto \sum_{\rho \in \Sigma(1)} \lfloor -\beta^* m(u_\rho) \rfloor D_\rho \]

To translate between Bondal's convention and Hanlon-Hicks-Lazarev's convention, one must multiply by $-1: \R^n \rightarrow \R^n$. With Hanlon-Hicks-Lazarev's convention, we have monomial labeling \[ p \mapsto \frac{ x^{\lfloor p \rfloor }  }{y^{ \lceil{p} \rceil}}. \]

When $p$ has integer coordinates (i.e., in the ``unimodular" setting of \cite{bayer-popescu-sturmfels}) then the Laurent monomial labeling from \cite{anderson2023resolutiondiagonalsmoothtoric} and \cite{hanlon2023resolutions} agree. When $X_\Sigma$ is smooth and non-unimodular, the resolutions will be different.

In \cite{bondalOberwolfach}, Bondal gave a criterion for the map $\Phi$ above which is used in the construction of the Hanlon-Hicks-Lazarev resolution of the diagonal to yield a full strong exceptional collection of line bundles on the smooth projective toric variety $X$ of dimension $n$: for each irreducible toric curve $C$, denote by $(D_1, \dots, D_{n-1})$ the irreducible toric divisors that contain $C$ and by $(a_1, \dots, a_{n-1})$ the corresponding intersection numbers with $C$. Bondal's numerical criterion is that all $a_i \geq -1$, with $a_i=-1$ appearing no more than once, which Bondal states holds for 16 out of 18 smooth projective toric Fano 3-folds.

\subsection{Relation between Bondal's numerical criterion and Bayer-Popescu-Sturmfels' condition of unimodularity}

\begin{defn} Let $X$ be a normal toric variety with fan $\Sigma$, and let $B = \left[ \begin{matrix} \leftarrow \rho_1 \rightarrow \\ \leftarrow\rho_2 \rightarrow \\ \vdots \\ \leftarrow\rho_n\rightarrow  \end{matrix} \right]$ be the matrix formed with rows given by the primitive ray generators $\{ \rho_i \text{ }|\text{ } 1\leq i \leq n \} $ where $|\Sigma(1)|=n$. Then $X$ is \textbf{unimodular} in the convention of \cite{bayer-popescu-sturmfels} if the matrix $B$ has linearly independent columns, and all maximal minors of $B$ lie in $\{0, \pm 1\}$. \end{defn} 

Throughout this paper, we use the condition ``unimodular" in this sense of \cite{bayer-popescu-sturmfels}. As Bayer-Popescu-Sturmfels note \cite[Section 6]{bayer-popescu-sturmfels}, the collection of line bundles $ \mathcal{E} = \{E_i\}$ appearing on one side of their resolution of $\mathcal{O}_\Delta$ need not yield a full strong exceptional collection, since for $X$ a unimodular toric variety, it is possible that $\mathcal{G} = \bigoplus_{i\in I}E_i$ can have $\mathrm{Ext}^i(\mathcal{G}, \mathcal{G}) \neq 0$ for some $i>0$. Note that the Hanlon-Hicks-Lazarev resolution of the diagonal yields an isomorphic resolution of $\mathcal{O}_\Delta$ to that of Bayer-Popescu-Sturmfels when $X$ is unimodular. 

 We verify that 16 out of 18 smooth projective toric Fano 3-folds are unimodular, and that this coincides with the success of the Hanlon-Hicks-Lazarev resolution of the diagonal to yield a full strong exceptional collection of line bundles. Furthermore, we also find that 96 out of 124 smooth projective toric Fano 4-folds are unimodular, and that 554 out of 866 smooth projective toric Fano 5-folds are unimodular. Macaulay2 code for checking the dimension 3-5 cases is given in the file in  ``UnimodCheckDim3To5.txt" in \eqref{link: githublink}, with output given by the length of AList describing the number of smooth projective toric Fano 5-folds which are non-unimodular. In agreement with \cite[Section 6]{bayer-popescu-sturmfels}, we find that the number of smooth projective toric Fano d-folds for which the Hanlon-Hicks-Lazarev resolution of the diagonal yields a full strong exceptional collection becomes strictly less than the number of unimodular smooth projective toric Fano d-folds for d=4,5. Under the assumption that Bondal's numerical criterion coincides with the success of the Hanlon-Hicks-Lazarev resolution of the diagonal to yield a full strong exceptional collection for smooth projective toric Fano varieties, this also implies that Bondal's numerical criterion is more restrictive than the condition of being unimodular for $X$ a smooth projective toric Fano variety in the sense that it holds for a lower proportion of all smooth projective toric Fano varieties of dimension $d$.

\section{Data collection}
Again, our data for checking whether the Hanlon-Hicks-Lazarev resolution of the diagonal yields a full strong exceptional collection on each of the 124 smooth projective toric Fano 5-folds is available at \newline \url{https://github.com/reggiea91/ADJOINT2025_Fivefolds.git}. For the sake of exposition, we include a negative and positive example, respectively, from our data collection procedure. 

\subsection{ Negative example in dimension 5 }

Let $X$ be the smooth projective toric Fano variety given by index 200 in the Macaulay2 database, accessed via \begin{verbatim} X = smoothFanoToricVariety(5,200)\end{verbatim} so that the complete fan $\Sigma$ has rays 

\[ \{(-1,0,0,0,0),(0,-1,0,0,0),(0,0,-1,0,0),(0,0,1,1,0),(0,0,0,-1,0),(0,0,0,0,1),(0,0,0,0,-1),(1,1,0,0,-2)\}\]

and presentation of the class group

\[  \Z^8 \stackrel{ \left( \begin{matrix} 
0&0&1&1&1&0&0&0\\
0&0&0&0&0&1&1&0\\
1&1&0&0&0&2&0&1\end{matrix} \right) }{\longrightarrow} \Z^3. \]
The Hanlon-Hicks-Lazarev resolution yields the free ranks (written as an ungraded resolution of S-modules, for S the
homogeneous coordinate ring of $Y \cong X \times X)$: 

\[0 \rightarrow S^{8} \rightarrow S^{34} \rightarrow S^{55} \rightarrow S^{42} \rightarrow S^{15} \rightarrow S^2 \rightarrow 0. \]

The Hanlon-Hicks-Lazarev resolution yields the collection of line bundles $\{\mathcal{O}(a_1, a_2, a_3)\}$ for $(a_1, a_2, a_3)$ appearing in

\begin{align*} 
\mathcal{E}  = & \{  (0, 0, 0), (0, -1, -1), (0, 0, -1), (-1, 0, 0), (-1, -1, -1), (0, -1, -2), (0, -1, -3), (-2, 0, 0), (-2, -1, -1), \dots \\ &\dots (0, 0, -2), (-1, 0, -1), (-1, -1, -2), (-1, 0, -2), (-2, -1, -2), (-1, -1, -3), (-2, 0, -1), (0, -1, -4), \dots \\ & \dots (-2, -1, -3), (-2, 0, -2), (-1, -1, -4), (-2, -1, -4) \}  \end{align*} 

on the left-hand side of the box product. Here, there does not exist an ordering for which $\mathcal{E}$ is exceptional since, for instance, the pair of line bundles $\mathcal{O}_X(0,-1,-1)$ and $\mathcal{O}_X(0,-1,-4)$ have non-zero homomorphisms in both directions. That is,

\[\Hom^\bullet_{D^b(X)}( \mathcal{O}_X(0,-1,-1), \mathcal{O}_X(0,-1,-4)) \cong H^*(X, \mathcal{O}(0,0,3)) \cong  \begin{cases} \C^{10} & \text{ in degree }0,\\ 0 & \text{ else. } \end{cases}    \] 

and

\[ \Hom^\bullet_{D^b(X)}( \mathcal{O}_X(0,-1,-4)), \mathcal{O}_X(0,-1,-1)) \cong H^*(X, \mathcal{O}(0,0,-3)) \cong  \begin{cases} \C & \text{ in degree } 2, \\ 0 & \text{ else}.        \end{cases}  \]

Here, we see that there does \textbf{not} exist an ordering of $\mathcal{E}$ for which we get an exceptional collection of line bundles. \\


\subsection{Positive example in dimension 5 }

Let $X$ denote the smooth projective toric variety given by index 851, so that 

\begin{verbatim}
    X = smoothFanoToricVariety(5,851)
\end{verbatim}

The complete fan $\Sigma$ has rays 

\[ \{(-1,0,0,0,0), (0,-1,0,0,0), (0,0,1,0,0), (0,0,-1,0,0), (0,0,0,-1,0), (0,0,0,1,1), (0,0,0,0,-1),(1,1,0,0,-1) \} \]

and we use as presentation of $\mathrm{Cl}(X)$:

\[ \Z^8 \stackrel{ \left( \begin{matrix}  0 &0 &1 &1 &0 &0 &0 &0 \\ 0 &0 &0 &0 &1 &1 &1 &0 \\ 1 &1 &0 &0 &1 &1 &0 &1   \end{matrix} \right) }{\longrightarrow} \Z^3. \]

The Hanlon-Hicks-Lazarev resolution yields the free ranks (written as an ungraded resolution of S-modules, for $S$ the homogeneous coordinate ring of $Y \cong X \times X)$:

\[ 0 \rightarrow S^6 \rightarrow S^{26} \rightarrow S^{43} \rightarrow S^{33} \rightarrow S^{11} \rightarrow S^1 \rightarrow  0.\]

There are 18 line bundles appearing on left-hand side which give a full strong exceptional collection of line bundles as follows:

\begin{align*} \mathcal{E} =\{ &  (-1,-2,-4) ,(0,-2,-4),(-1,-2,-3),(-1,-1,-3),(0,-2,-3),(-1,-2,-2),(-1,-2,-2),(-1,-1,-2)\dots \\ &\dots(0,-1,-3),(-1,0,-2),(-1,-1,-1),(-1,0,-1),(0,-2,-2),(0,-1,-2),(0,0,-2),(-1,0,0),(0,-1,-1),\dots \\ &\dots(0,0,-1),(0,0,0) \}. 
  \end{align*}

Here $\Hom_{\mathrm{D}^b_{Coh}(X)}(E_i, E_j)$ is concentrated in degree 0 for all $i$ and $j$, with the rank of $\Hom^0(E_j, E_i)$ given by the $(i,j)$ entry of the following matrix: 

\begin{align*}
\left( \begin{matrix} 1 &0 &0 &0 &0 &0 &0 &0 &0 &0 &0 &0 &0 &0 &0 &0 &0 &0 \\
3 &1 &0 &0 &0 &0 &0 &0 &0 &0 &0 &0 &0 &0 &0 &0 &0 &0 \\
5 &1 &1 &0 &0 &0 &0 &0 &0 &0 &0 &0 &0 &0 &0 &0 &0 &0 \\
2 &0 &0 &1 &0 &0 &0 &0 &0 &0 &0 &0 &0 &0 &0 &0 &0 &0 \\
6 &3 &0 &0 &1 &0 &0 &0 &0 &0 &0 &0 &0 &0 &0 &0 &0 &0\\
12 & 5 &3 &0 &1 &1 &0 &0 &0 &0 & 0&0 &0 &0 &0 &0 &0 &0 \\
15 &5 &5&0 &1 &1& 1& 0& 0& 0 & 0& 0 &0 &0 &0 &0 &0 &0\\
6 &2 &0 &3 &0 &0 &0 &1 &0 &0 &0 &0 &0 &0 &0 &0 &0 &0\\
10 &2 &2 &5 &0 &0 &0 &1 &1 &0 &0 &0 &0 &0 &0 &0 &0 &0\\
12 &6 &0 &6 &2 &0 &0 &3 &0 &1 &0 &0 &0 &0 &0 &0 &0 &0 \\
22 &12 &6 &0 &5 &3 &0 &0 &0 &0 &1 &0 &0 &0 &0 &0 &0 &0 \\
24 &10 &6 &12 &2 &2 &0 &5 &3 &1 &0 &1 &0 &0 &0 &0 &0 &0\\
30 &10& 10& 15& 2 &2 &2 &5 &5 &1 &0 &1 &1 &0 &0 &0 &0 &0\\
31 &15 &12& 0& 5& 5 &3 &0 &0 &0 &1 &0 &0 &1 &0 &0 &0 &0\\
44 &24& 12 &22& 10& 6 &0 &12 &6 &5 &2 &3 &0 &0 &1 &0 &0 &0\\
62 &30& 24& 31& 10& 10& 6 &15& 12& 5& 2& 5& 3 &2 &1 &1 &0 &0\\
53 &31 &22 &0 &15& 12 &6 &0 &0 &0 &5 &0 &0 &3 &0 &0 &1 &0\\
106 &62& 44 &53& 30 &24& 12 &31& 22& 15& 10& 12& 6& 6& 5& 3& 2& 1 
\end{matrix} \right) 
\end{align*} 

Exceptionality of $\mathcal{E}$ is now indicated by the fact that this matrix is lower-diagonal.

\section{Future directions} 
In dimension 1, the Hanlon-Hicks-Lazarev resolution yields a full strong exceptional collection of line bundles for the unique smooth projective toric Fano curve $\PP^1$. For $\PP^n$, the resulting collection of line bundles agrees with the collection due to Beilinson. In dimension 2, the resolution of the diagonal developed by the author coincides with the Hanlon-Hicks-Lazarev resolution of the diagonal for all 5 smooth projective toric Fano surfaces, since all are unimodular in the sense of Bayer-Popescu-Sturmfels. In this case, either resolution of the diagonal yields a full strong exceptional collection for all 5 smooth projective toric Fano surfaces \cite{anderson2024exceptionalcollectionslinebundles}. In dimension 3, 2 out of 18 smooth projective toric Fano 3-folds provably do not admit a full strong exceptional collection of line bundles from the Hanlon-Hicks-Lazarev resolution of the diagonal due to failing Bondal's numerical condition \cite{bondalOberwolfach}. In dimension 4, the Hanlon-Hicks-Lazarev resolution of the diagonal yields a full strong exceptional collection for 72 out of 124 smooth projective toric Fano fourfolds \cite{ramirez2025exceptionalcollectionstoricfano}. We have found that for 300 out of 866 smooth projective toric Fano fivefolds, the Hanlon-Hicks-Lazarev resolution of the diagonal yields a full strong exceptional collection of line bundles. This leads to the following conjecture.

\begin{conjecture}

Let $f(d): \Z_{>0} \rightarrow [0,1]$ denote the proportion of smooth projective toric Fano $d$-folds for which the Hanlon-Hicks-Lazarev resolution of the diagonal yields a full strong exceptional collection of line bundles. Then we have 

\[ \frac{1}{1}= f(1)=f(2) = \frac{5}{5},\hspace{.2cm}  f(3) = \frac{16}{18}, \hspace{.2cm}  f(4) = \frac{72}{124}, \hspace{.2cm} f(5) =  \frac{300}{866}. \]

\flushleft 
We conjecture that $f(d)$ is strictly monotone decreasing for $d>1$, and that $\lim_{d\rightarrow \infty} f(d) = 0$. \end{conjecture}

\section{Acknowledgements} We are grateful to Jay Yang for providing Macaulay2 code reproducing the Hanlon-Hicks-Lazarev resolution, with contributions to the code from Mahrud Sayrafi. 
This material is based upon work supported by the Sloan foundation under Grant G-2025-25245 while the author participated in the Self-ADJOINT program hosted by the Simons Laufer Mathematical Sciences Institute in Berkeley, California in July 2025. The author also thanks Candace Bethea and Alicia Lemarche for their contributions to this project.  

\section{Declaration of Interest}
The author reports there are no competing interests to declare.

\printbibliography

\end{document}